\input amstex.tex
\input amsppt.sty
\documentstyle{amsppt}
\NoBlackBoxes
\magnification=\magstep1
\voffset=-1.3cm
\advance\hsize0.5cm

\def\delet{\mathaccent"7017 }

\def\Ker{\mathop{\fam0 Ker}}

\def\rmop#1{\expandafter\def\csname#1\endcsname{\operatorname{#1}}}
\def\R{\Bbb R} \def\Z{\Bbb Z}  \def\Sy{\Bbb S} \def\P{\Bbb P}

\def\t{\widetilde}

\topmatter
\title On the Browder-Levine-Novikov embedding theorems \endtitle
\author M. Cencelj, D. Repov\v s and A. Skopenkov \endauthor
\address M. Cencelj and D. Repov\v s: Faculty of Education and Faculty of Mathematics and Physics, University of Ljubljana \& Institute for Mathematics, Physics and Mechanics, Slovenia.
\newline
E-mail: matija.cencelj\@fmf.uni-lj.si, dusan.repovs\@fmf.uni-lj.si
\newline
Web: https://www.fmf.uni-lj.si/en/directory/494/cencelj-matija/   
\newline
https://www.fmf.uni-lj.si/en/directory/206/repovs-dusan/
\endaddress
\address
A. Skopenkov: Moscow Institute of Physics and Technology, and Independent University of Moscow, Russia.
\newline
Email: skopenko\@mccme.ru. Web: https://users.mccme.ru/skopenko/ \endaddress
\thanks This is an updated version of a survey published in Proc. Steklov Math. Inst. 247 (2004).
We put it on arXiv because the survey attracted some attention.
We would like to acknowledge P. Lambrechts and Yu. P. Solovyov for useful discussions and remarks.
\newline
{\it 2010 Math. Subj. Classification:} Primary: 57R40; Secondary:  57R42, 57R52, 57R65, 57Q35, 57Q37.
\endthanks

\abstract
In this survey we present applications of the ideas of complement and neighborhood in the theory embeddings of manifolds into Euclidean space (in codimension at least three).
We describe how the combination of these ideas gives a reduction of embeddability and isotopy problems to algebraic problems.
We present a more clarified exposition of the Browder-Levine theorem on realization of normal systems.
Most of the survey is accessible to non-specialists in the theory of embeddings.
\endabstract


\endtopmatter

\document

\head Introduction \endhead

Three important classical problems in topology are the following, cf. [Ze93, p. 3].

(1) {\it The Manifold Problem:} find conditions under which two spaces are homeomorphic (and also describe the homeomorphism classes of manifolds from a given class);

(2) {\it The Embedding Problem:} find conditions under which a space embeds into $S^m$ for a given $m$; and

(3) {\it The Knotting Problem:} find conditions under which two embeddings
are isotopic (and also describe the isotopy classes of embeddings $N\to S^m$).

This survey is concerned with the second and the third problem.
We show how the ideas of complement and neighborhood can be used to study embeddings.
Then we combine these two ideas and formulate the Browder-Levine-Novikov Theorems 6 and 8, and the Browder-Wall Conjecture 9.
Such an exposition of these results has apparently not been published yet.

In this survey we present classical theorems that give readily calculable results [Sk16, Remark 1.2], at least for particular cases.
So we do not mention other results like Browder-Haefliger-Casson-Sullivan-Wall theorem stating that {\it a homotopy equivalence between a closed $k$-manifold and an $m$-manifold with boundary is homotopic to an embedding (in the smooth category for $2m\ge3k+3$, in the Pl category for $m\ge k+3$).}
The Browder-Levine-Novikov Theorems 6 and 8 give readily calculable results [Sk16, Remark 1.2], although only for simpler manifolds like homotopy spheres.
Also the classification theorem for links [Ha66] can be proved analogously to the Browder-Levine Theorem 6 [Ha86].
In general, `these theorems reduce geometric problems to algebraic problems which are even harder to solve' [Wa70].
However, where the Browder-Levine-Novikov Theorems fail to give readily calculable results, the modified Kreck surgery could be applied, see [Sk05, Sk06', CS08, CS16I, CS16II].
For surveys on different approaches to the Embedding and the Knotting Problems see [RS99, Sk06, Sk16].

Denote by CAT the smooth (DIFF) or piecewise-linear (PL) category.
We omit CAT if a definition or a statement holds in both categories.
A polyhedron $N$ is called {\bf PL embeddable} into $ S^m$, if there is a PL injective map $f:N\to S^m$.
A smooth manifold $N$ is called {\bf smoothly embeddable} into $S^m$,
if there is a smooth injective map $f:N\to S^m$ such that $df(x)$ is nondegenerate for each $x\in N$.
Such a map $f$ is called an {\bf embedding} of $N$ into $ S^m$ (in the corresponding category).

Two embeddings $f,f':N\to S^m$ are said to be {\bf (ambient) isotopic},
if there exists a homeomorphism onto $F: S^m\times I\to S^m\times I$ such that

(i) $F(y,0)=(y,0),\quad\text{for each}\quad y\in S^m$;

(ii) $F(f(x),1)=f'(x),\quad\text{for each}\quad x\in N$; \quad and

(iii) $F(S^m\times\{t\})= S^m\times\{t\},\quad\text{for each}\quad t\in I.$

This homeomorphism $F$ is called an {\it (ambient) isotopy}.
An (ambient) isotopy is also a homotopy $S^m\times I\to S^m$ or a family of
maps $F_t:S^m\to S^m$, generated by the map $F$ in the obvious way.
See [Is] for more information.


\head The idea of complement and obstructions to embeddability \endhead

In some cases one can obtain obstructions to embeddability using the
following idea which can be traced back to early works of J.W. Alexander (around 1910).
Considering the complement $S^m-N$ of $N\subset S^m$,
one can deduce necessary conditions for $N$ itself.

Let us illustrate this idea by embeddability of graphs in the plane.
Suppose that a connected graph $N$ embeds into the plane.
By the Euler formula, $V-E+F=2$.
Since every face has at least three edges on its boundary, we have $3F\le2E$.
Therefore $E\le3V-6$.
This implies that the graph $K_5$ is not planar.
Analogously, one can prove that the graph $K_{3,3}$ is not planar.


In general, the Euler formula is replaced by its generalization, that is, by the Alexander duality.
For example, if $N\subset S^m$, then the Betti numbers satisfy $b^m(N)=b^{-1}(S^m-N)=0$.
This gives a necessary condition for $N$ to be embeddable into $S^m$.
By developing this idea, W. Hantzsche proved that $\R P^3$ and the 3-dimensional lens spaces do not embed into $S^4$, and that $\R P^{4k-1}$ does not embed into $S^{4k}$ (this follows by Theorem 1.b below).

\proclaim{Theorem 1} [Ha37, Theorems 1 and 3] 
(a) If $N$ is a closed $2l$-submanifold of $S^{2l+1}$, then $N$ is orientable and the Euler characteristic of $N$ is even (i.e. the dimension of the free part of $H_l(N;\Z)$ or the $l$-th Betti number of $N$ or $\dim H_l(N;\Z_2)$ is even).

(b) If $N$ is a closed $(2l+1)$-submanifold of $S^{2l+2}$, then $N$ is orientable
and the torsion part of $H_l(N;\Z)$ is isomorphic to $G\oplus G$ for some finite abelian group $G$.
\endproclaim

{\it Proof of (b).} Denote by $A$ and $B$ the closures of the complement in $S^{2l+2}$ of $S^{2l+2}-N$.
In this proof we omit $\Z$-coefficients.
Using the Mayer-Vietoris sequence for $S^{2l+2}=A\cup B$, $N=A\cap B$, we see that the sum of inclusions induces an isomorphism $H_l(N)\cong H_l(A)\oplus H_l(B)$.
For the torsion subgroups we have
$$TH_l(B)\overset{L}\to\cong TH_{l+1}(B,\partial B)\overset{E}\to\cong TH_{l+1}(S^{2l+2},A)\overset{\partial}\to
\cong TH_l(A).$$
Here $L$ is the Lefschetz duality isomorphism, $E$ is the excision isomorphism and $\partial$ is the isomorphism from the exact sequence of the pair $(S^{2l+2},A)$.
Hence $TH_l(N)\cong G\oplus G$ for $G=TH_l(A)$.
\qed



\smallskip
Part (a) is better known and holds under the weaker assumption that $N$ is the boundary of a manifold.  

Further, by studying duality between homology {\it rings} of a manifold $N\subset S^m$ and the complement $S^m-N$ H. Hopf proved that $\R P^{m-1}$ does not embed into $S^m$ (for the remaining case $m=4k+2$)\footnote{It is not explained in [SE62] how this non-embeddability [SE62, III, Corollary 2.2] follows from [SE62, III, Theorem 2.1] for the non-trivial case when $n=m$ even. In the notation of [SE62] it is unclear why we cannot have $i^*H^q(A)=H^q(M)$ and $j^*H^q(B)=0$ for $0<q<n/2$, $i^*H^q(A)=0$ and $j^*H^q(B)=H^q(M)$ for $n/2\le q<n-1$.
The first statement of [SE62, III, Theorem 2.1] is the (dual) isomorphism $H_l(N)\cong H_l(A)\oplus H_l(B)$ from the above proof of Theorem 1.b.
Theorem 1.b gives non-embeddability of $\R P^{m-1}$ into $S^m$ for $m=4k$ but not for $m=4k+2$. 
Theorem 2.b below is a simplified version  of [SE62, III, Theorem 2.1] perhaps sufficient for applications.}
[Ho40], [Sk20, solution of Problem 11.1.d].
See [ARS], [Sk20, Problem 11.1.e] for a simple development of this idea proving that
$\R P^3\times \R P^3$ does not embed into $S^7$.
A simpler version of this idea implies that $\R P^3\#\R P^3$ does not embed into $S^4$ (this follows from Theorem 2.b below but not from Theorem 1.b).

\proclaim{Theorem 2}
(a) If $N$ is a $2l$-submanifold of $S^{2l+1}$, then $H_l(N;\Z_2)$ is the sum $A\oplus B$ for some isomorphic subgroups $A,B$ such that the modulo 2 intersection form of $N$ [IF] vanishes on each them
(and so the form is non-degenerate on $A\times B$).

(b) If $N$ is a $(2l+1)$-submanifold $N$ of $S^{2l+2}$, then the torsion part $T$ of $H_l(N;\Z)$ is the sum $A\oplus B$ for some isomorphic subgroups $A,B$ of $T$ such that the linking form of $N$ [LF] vanishes on each of them (and so the form is non-degenerate on $A\times B$). [GS99, Exercise 4.5.12.d]
\endproclaim

Part (a) is a $\Z_2$-analogue of the well-known $\Z$-version for $l$ even.
The version is equivalent to the statement that the signature of a boundary (in particular, of a closed codimension 1 submanifold of Euclidean space) is zero.

R. Thom obtained some conditions on the cohomology ring of a closed $(m-1)$-manifold,
which are necessary for embeddabilty into $S^m$ [Th51].


\head The idea of complement and obstructions to isotopy \endhead

The idea of the complement is even better applicable to study the Knotting Problem.
Indeed, if $f,g:N\to S^m$ are isotopic embeddings, then $S^m-f(N)\cong S^m-g(N)$.
Therefore any invariant of the space $S^m-f(N)$ is an isotopy invariant.

This idea was first applied by Alexander in 1910's to knots in 3-space.
Let us illustrate his idea by the proof that the trefoil knot $f:S^1\to S^3$ is not isotopic to the trivial knot.
Using the Van Kampen theorem on the fundamental group of the union, we obtain that
$\pi:=\pi_1(S^3-f(S^1))=\left<x,y\ |\ xyx=yxy\right>$.
Evidently, $\pi_1(S^3-S^1)=\Z$.
The first idea to distinguish between $\Z$ and $\pi$ is to compare $\Z$ and
the {\it abelianization} $\pi/[\pi,\pi]$ of $\pi$.
But it turns out that $\pi/[\pi,\pi]=\Z$.
Perhaps Alexander, trying so to distinguish knots, observed that
$\pi/[\pi,\pi]=\Z$ for fundamental group $\pi$ of the complement to {\it any}
knot, which lead him to discovery of his duality theorem.
In order to distinguish between the trefoil knot and the trivial knot, one can
construct a non-trivial homomorphism $\pi\to S_3$, defined by $x\to(12)$, $y\to(23)$.
Hence $\pi$ is not abelian and not isomorphic to $\Z$.

Exposition of the subsequent development of the complement idea in the theory of knots $S^1\subset S^3$ (or more generally, $S^n\subset S^{n+2}$) is beyond the purposes of this survey (for a user's guide see [Sk20u]).
We only formulate sufficient conditions for the completeness of the complement
invariant (without proofs).

\proclaim{Theorem 3}  For any $n\ne2$, a smooth embedding $f:S^n\to S^{n+2}$ is smoothly isotopic to the standard embedding if and only if $S^{n+2}-f(S^n)\simeq S^1$ [Le65].
\endproclaim

{\bf Remark.} (a) By [Pa57] $(S^3-f(S^1))\simeq S^1$ is equivalent to $\pi_1(S^3-f(S^1))\cong\Z$.

(b) An analogue of Theorem 3 also holds in PL locally flat and TOP locally flat
categories for $n\ne2,3$ [Pa57, Le65, St63], see also [Gl62].
Recall that a (PL or TOP) embedding $S^n\subset S^m$ is said to be (PL or
TOP) {\it locally flat}, if every point of $S^n$ has a neighborhood $U$ in
$S^m$ such that $(U,U\cap S^n)$ is (PL or TOP)
homeomorphic to $(B^n\times B^{m-n},B^n\times0)$.

The local flatness assumption in TOP category is necessary in order to rule
out {\it wild} embeddings, which were first constructed by Antoine in 1920
(the Antoine necklace) and Alexander in 1923 (the Alexander Horned Sphere),
again using the complement idea [Al24'], see also 
[DV09].

(c) A motivation for the construction of the Alexander Horned sphere is study
of the TOP knots in codimension 1, i.e. of embeddings $S^n\to S^{n+1}$.
A short history of this problem is as follows.
The well-known Jordan theorem, first proved by Brouwer, states that every
$S^n$, contained in $S^{n+1}$, separates $S^{n+1}$ into two components.
It is easy to prove the `analogue' of Theorem 3: $S^n\subset S^{n+1}$ is
unknotted if and only if the closures of these components are balls.
In 1912 Sch\"onflies proved that every $S^1\subset S^2$ is unknotted.
Thus, unknottedness $S^n\subset S^{n+1}$ is called the `Sch\"onflies theorem' or the `problem'.
In 1921 Alexander announced that he has proved the Sch\"onflies theorem for arbitrary $n$.
However, in 1923 he found a counterexample --- the celebrated Horned Sphere [Al24'].
Nevertheless, he proved the PL Sch\"onflies theorem for $n=2$ [Al24].
For $n\ge3$ the PL Sch\"onfliess conjecture is a famous difficult unsolved problem [RS72].

\proclaim{Theorem 4} Any two smooth embeddings $S^n\to S^{n+1}$ are smoothly isotopic for $n\ne3$ [Sm61, Ba65].
\endproclaim


After the appearance of the Alexander Horned Sphere the unknottedness of a PL locally flat or TOP locally flat embedding $S^n\to S^{n+1}$ was also called the {\it Sch\"onflies conjecture}.
It was proved only in 1960s
(in the PL case only for $n\ne3$) [Ma59, Br60, Mo60, RS72].
Note that Brown's elegant short proof [Br60] of the TOP locally flat analogue of Theorem 4 gave rise to the theory of `cellular sets', which became an important branch of geometric topology.

\head The complement invariant \endhead


{\it In the rest of this paper, $N$ is a closed smooth $n$-manifold,
$f:N\to S^m=S^{n+k}$ is a smooth embedding, and $i:S^m\to S^{m+1}$ is the standard inclusion.}

Define the {\bf complement} to be
$$C(f)=S^m-f(N)\simeq S^m-Of(N),$$
where $Of(N)$ is an open tubular neighborhood in $S^m$ of $f(N)$.

The topological and the homotopy types of $C(f)$ are invariants of $f$.
For $m-n\ge3$, the complement $C(f)$ is simply-connected.
Hence by the Alexander duality, its homology groups do not depend on $f$.
Therefore the invariant $C(f)$ is relatively weak for $m-n\ge3$.

\smallskip
{\bf Example.} (a) Assume that $N=S^{n_1}\sqcup\dots\sqcup S^{n_r}$ and $f:N\to S^m$ is any embedding.
For $m-n_i\ge3$, the complement $C(f)$ is simply-connected.
Hence by the Alexander duality and Whitehead Theorem $C(f)\simeq S^{m-n_1-1}\vee\dots\vee S^{m-n_r-1}$.

(b) {\it The complement of $N=S^p\times S^q$ may depend on the embedding into $S^m$ even for $m-p-q\ge3$.}
Indeed, consider the Hopf fibration $S^3\to S^7\overset{h}\to\to S^4$.
Take the standard embedding $S^2\subset S^4$.
Its complement has the homotopy type of $S^1$.
We have $h^{-1}(S^1)\cong S^1\times S^3\subset S^7$, and the complement
of this embedding $f$ is
$$C(f)\simeq h^{-1}(S^2)\cong S^2\times S^3\not\simeq S^2\vee S^3\vee
S^5\simeq C(f_0),$$
where $f_0:S^1\times S^3\to S^7$ is the standard embedding.
Analogously, one can construct two embeddings $S^3\times S^7\subset S^{15}$
whose complements are not homotopy equivalent.
This example is due to P. Lambrechts (personal communication).

\smallskip
By general position, $C(f)$ does not depend on $f$ for $m\ge2n+2$.
This space is denoted by $C_m(N)$.
We have
$$C(i\circ f)\simeq\Sigma C(f),\quad\text{so}\quad
\Sigma^{M-m}C(f)\simeq C_M(f)\quad\text{for}\quad M\ge2n+2.$$
The {\it necessary complement condition} for embeddability of $N$ into $S^m$
is the $(M-m)$-desuspendability of $C_M(N)$, i.e. the existence of a space
$C$ such that $\Sigma^{M-m}C\simeq C_M(f)$.
If this condition holds for $M=2n+2$, then it automatically holds for each
$M\ge2n+2$.

\head The idea of neighborhood \endhead

In order to obtain necessary conditions for embeddability into $S^m$ one can assume that $N\subset S^m$ and consider relations between $N$ and its neighborhood in $S^m$.
This method seems to have been first introduced by H. Whitney in the smooth category.
He proved many non-trivial results, like non-embeddability of $\R P^n$ into $\R^{2n-1}$ for $n$ a power of 2
(see a user-oriented exposition in [Sk20, \S12]).
For this, he created a theory of bundles and introduced the so-called Stiefel--Whitney classes
$w_k(N)\in H^k(N;\Z_2)$ and the dual Stiefel--Whitney classes $\overline w_k(N)\in H^k(N;\Z_2)$ [SW].
These classes have played an important role in topology and differential geometry.

A smooth map $g:N\to S^m$ is a {\it smooth immersion} if $dg(x)$ is nondegenerate for each $x\in N$.
Two smooth immersions $g,g':N\to S^m$ are {\it regularly homotopic},
if there exists a smooth immersion $G:N\times I\to S^m\times I$ such that

(i) $G(x,0)=(g(x),0)$, $G(x,1)=(g'(x),1)$ for each $x\in N$; and

(ii) $G(N\times\{t\})\subset S^m\times\{t\}$ for each $t\in I$.

The normal bundle $\nu(g)$ of $g$ up to equivalence is a regular homotopy invariant of $g$,
and an isotopy invariant of $g$, if $g$ is an embedding.
This isotopy invariant is not very strong, because e.g. normal bundles of different embeddings are stably equivalent.

For an embedding $f:N\to S^{n+k}$ we may assume that the space of $\nu(f)$ is $Of(N)$ and the zero section is $f$.
By general position for $k\ge n+2$  bundles $\nu(f)$ are isomorphic for different $f$.
Each such bundle is called {\it the stable normal bundle} $\nu_k(N)$ of $N$.
We have 
$$\nu(i\circ f)=\nu(f)\oplus1,\quad\text{so}\quad
\nu(f)\oplus(K-k)=\nu_K(f)\quad\text{for}\quad K\ge n+2.$$
The {\it necessary normal bundle condition} for immersability of $N$ to $S^{n+k}$ is the existence of a $k$-bundle $\nu$ over $N$ that is stably equivalent to $\nu_K(f)$.
Since $\nu(f)\oplus\tau(N)=n+k$, it follows that the normal bundle condition is equivalent to the stable triviality of $\nu\oplus\tau(N)$, and then to the triviality of $\nu\oplus\tau(N)$ (by general position, because the dimension of $\nu\oplus\tau(N)$ is greater than $n$).


\proclaim{Smale-Hirsch Theorem 5} [Hi59,  Im]
If there exists a $k$-bundle $\nu$ over $N$ such that
$\nu\oplus(K-k)\cong\nu_K(N)$, then $N$ immerses in $\R^{n+k}$.

If $f,g:N\to\R^{n+k}$ are immersions such that $df,dg:TN\to\R^{n+k}$ are homotopic
via linear monomorphisms, then $f$ and $g$ are regular homotopic.
\endproclaim

The notion of Stiefel--Whitney classes was generalized by L.S. Pontryagin who introduced
Pontrjagin classes $p_{4k}(N)\in H^{4k}(N;\Z)$ and their duals $\overline p_{4k}(N)\in H^{4k}(N;\Z)$ for oriented $N$.
M. Atiyah and F. Hirzebruch studied normal bundles as elements of $K(N)$ and obtained interesting non-embedding results [FF89, \S41.6].

\head Normal systems and the Browder-Levine embedding theorem \endhead

Combining the complement and the neighborhood ideas, Levine, Novikov and Browder obtained
new necessary conditions for embeddability of manifolds (and for isotopy of embeddings).
The proofs of sufficiency of these conditions are one of the most important
applications of {\it surgery} to topology of manifolds (see also [Ha66, Ha68, Gi71, Sk05, \S3]).
For a closely related concept of thickenability and thickenings see [Wa66, LS69, RBS].

Denote by $a(f)$ the homotopy class of the inclusion $\partial C(f)\subset C(f)$.
This class is an element of a set depending on $f$, so $a(f)$ alone is not an invariant of $f$.
The corresponding (`attaching') invariant of $f$ is defined as follows.

Regard $\nu(f)$ as the
normal $S^{k-1}$-bundle $\partial C(f)\to N$.

The {\bf normal system of $f$} is the triple $\Sy(f)=(\nu(f),C(f),a(f))$.
A {\bf normal system of codimension $k$ on a manifold $N$} is the triple $\Sy=(\nu,C,a)$ of
an $S^{k-1}$-bundle $\nu:E\to N$, a space $C$ and a map $a:E\to C$.

Two normal systems $(\nu,C,a)$ and $(\nu_1,C_1,a_1)$ are called {\bf equivalent}, if
there is a bundle isomorphism $b:E\to E_1$ and a homotopy equivalence $r:C\to C_1$ such that
$r\circ p\simeq p_1\circ b$.

The normal systems of isotopic embeddings are equivalent.
So the equivalence class of the normal system is an isotopy invariant.


For $K>n+1$  every two embeddings $N\to S^{n+K}$ are isotopic,
so the normal systems $\Sy(f)$ are equivalent for different $f$.
Each such system is called the {\bf stable normal system} $\Sy_K(N)$ of $N$.

The {\it suspension} over a normal system $\Sy=(\nu,C,p)$ is the normal system
$\Sigma\Sy:=(\nu\oplus1,\Sigma C,a')$, where $a'$ is the suspension of $a$ on each fiber.
We have $\Sigma^{K-k}\Sy(f)=\Sy_K(N)$ for $K>n+1$.
{\it The Browder-Levine necessary condition} for embeddability of $N$ into $S^{n+k}$ is as follows.

{\bf BL(k)} {\it there is a normal system $\Sy=(\nu,C,a)$ of codimension $k$ on $N$ such that $\Sigma^{K-k}\Sy\sim\Sy_K(N)$.}

Clearly, the latter condition does not depend on $K>n+1$.

\proclaim{Browder-Levine Theorem 6}
Suppose that $K>n+1$, $\pi_1(N)=0$, $k\ge3$, BL(k) holds and $\pi_1(C)=0$.
Then there exists a smooth embedding $f:N\to S^{n+k}$ such that $\Sy(f)\sim\Sy$.
\endproclaim

Theorem 6 was proved in [Le65] for $N$ a homotopy sphere and in [Br68] for the general case.
Theorem 6 easily follows by induction from the following lemma on compression and desuspension.

\proclaim{Lemma 7}
Suppose that $F:N\to S^{n+k+1}$ is a smooth embedding, $\pi_1(N)=0$, $k\ge3$, $n+k\ge5$ and $\Sy=(\nu,C,a)$ is a normal system of codimension $k$ over $N$ such that $\Sigma\Sy\sim\Sy(F)$ and $\pi_1(C)=0$.
Then there exists a smooth embedding $f:N\to S^{n+k}$ such that $\Sy(f)\sim\Sy$.
\endproclaim

Proof of Lemma 7 is given at the end of the paper.
It uses the notion of normal invariant defined in the next section.

We conjecture that {\it if $n\ge 2k+4$ and $f,g:N\to S^{n+k}$ are smooth embeddings such that $\Sy(f)\sim\Sy(g)$, then $f$ and $g$ are smoothly isotopic.}
Cf. the Browder-Wall Conjecture 9.
These conjectures, even if they are correct, would hardly give readily calculable classification results, see the citation of [Wa70] in the introduction and [Sk05, the paragraph after proof of Lemma 1.3].
The required techniques were first introduced for the Manifold Problem.
They are extended to the Embedding Problem in the proof of Lemma 7, and could be extended to the Knotting Problem using the following well-known lemma (for a simple proof see [Sk05, Lemma 1.3]).

{\it For a closed connected manifold $N$ embeddings $f,f':N\to\R^m$ are
isotopic if and only if there is an orientation-preserving bundle isomorphism
$\varphi:\partial C(f)\to\partial C(f')$ which extends to an
orientation-preserving diffeomorphism $C(f)\to C(f')$.}

\head Normal invariants and the Browder-Novikov embedding theorem \endhead

Denote by $\alpha(f)$ the homotopy class of the composition
$$S^{n+k}\to S^{n+k}/C(f)\cong T\nu(f)\ \in\ \pi_{n+k}(T\nu(f))$$
of the quotient map and  canonical homeomorphism.
The class $\alpha(f)$ is an element of a group depending on $f$, so $\alpha(f)$ alone is not an invariant of $f$.
One introduces the corresponding invariant as follows (analogously to the previous section).

{\it The normal invariant of $f$} is the pair $(\nu(f),\alpha(f))$.
(Sometimes $\alpha(f)$ alone is called the normal invariant.)
{\it A normal invariant of codimension $k$ on a manifold $N$} is a pair $(\nu,\alpha)$ of
an $S^{k-1}$-bundle $\nu:E\to N$ and an element $\alpha\in\pi_{n+k}(T\nu)$.

Two normal invariants $(\nu,\alpha)$ and $(\nu_1,\alpha_1)$ are called {\it equivalent}, if
there is a bundle isomorphism $b:E\to E_1$ such that $\alpha_1=b_T\circ\alpha$, where $b_T$ is the map of Thom spaces corresponding to a bundle map $b$.

The normal invariants of isotopic embeddings are equivalent.
So the equivalence class of the normal invariant is an isotopy invariant.


For $K>n+1$ every two embeddings $N\to S^{n+K}$ are isotopic,
so their normal invariants are equivalent for different $f$.
Each such normal invariant is called {\it stable normal invariant} $(\nu_K(N),\alpha_K(N))$ of $N$.

We have
$$T(\nu\oplus1)\simeq\Sigma T(\nu)\quad\text{and}\quad
\alpha(i\circ f)=\Sigma\alpha(f).$$
The {\it suspension} over a normal invariant $(\nu,\alpha)$ is the normal invariant
$\Sigma(\nu,\alpha):=(\nu\oplus1,\Sigma\alpha)$.
We have $\Sigma^{K-k}(\nu(f),\alpha(f))=(\nu_K(N),\alpha_K(N))$ for $K>n+1$.
{\it The Browder-Novikov necessary condition} for embeddability of $N$ into $S^{n+k}$ is as follows.

{\bf BN(k)} {\it there is a normal invariant $(\nu,\alpha)$ of codimension $k$ on $N$
such that $\Sigma^{K-k}(\nu,\alpha)\sim (\nu_K(N),\alpha_K(N))$.}


Clearly, the latter condition does not depend on $K>n+1$.

\proclaim{Browder-Novikov Theorem 8}
Suppose that $K>n+1$, $k\ge2$, $\pi_1(N)=0$ and BN(k) holds.
Then there exists a smooth embedding $f:N\to S^{n+k+1}$ such that $(\nu(f),\alpha(f))\sim\Sigma(\nu,\alpha)$.
\endproclaim

Theorem 8 was proved in [No64] for $N$ a homotopy sphere, and in [Br68] for the general case (see also [Le63]).
Theorem 8 follows from the Browder-Levine Theorem 6 because BN(k)$\Rightarrow$BL(k+1).

\smallskip
{\it Proof that BN(k)$\Rightarrow$BL(k+1).}
Take $\nu$ and $\alpha$ from BN(k).
Define the following normal system from BL(k+1):
$$\Sy=(\nu\oplus1,C,a),\quad\text{where}
\quad C=T\nu\bigcup\limits_{\alpha:\partial D^{n+k+1}\to T\nu}D^{n+k+1}$$
$$\text{and $a$ is the composition}\quad S(\nu\oplus1)
\to S(\nu\oplus1)/\text{(the canonical section)}\cong T\nu\subset C.$$

Denote by $D\nu$ the space of the $D^k$-bundle associated to an $S^{k-1}$-bundle $\nu$.

\smallskip
{\it Proof that BL(k)$\Rightarrow$BN(k) for $N$ simply-connected and $k\ge3$.}
Take $\Sy=(\nu,C,a)$ from BL(k).
Since
$$\Sigma^{K-k}\Sy\sim\Sy(f),\quad\text{we have}
\quad D(\nu\oplus(K-k))\bigcup\limits_{a^{(K-k)}}\Sigma^{K-k}C\simeq S^{n+K}.$$
Since $N$ is simply-connected and $k\ge3$, both $D\nu\simeq N$ and $C$ are simply-connected.
Therefore using Mayer-Vietoris sequence, one can show that $S^{n+k}\simeq D\nu\cup_aC$.
The composition of this homotopy equivalence and collapsing of $C$ to a point represents the required $\alpha$.

\head Poincar\'e embeddings and the Browder-Wall embedding conjecture \endhead

Let us generalize the above results.
For an embedding $f:N\to M$ of a smooth manifold $N$ into a smooth manifold
$M$ define $C(f)$, $\nu(f)$ and $a(f)$ as in the above definition of the normal system of $f$.
Take the identity map $h(f):D\nu(f)\cup_{a(f)}C(f)\to M$.
The {\bf Poincar\'e embedding of embedding $f$} is the quadruple $\P(f)=(\nu(f),C(f),a(f),h(f))$.
A {\bf Poincar\'e embedding of a manifold $N$ into a manifold $M$} is a quadruple $\P=(\nu,C,a,h)$ consisting of
an $S^{k-1}$-bundle $\nu:E\to N$,
a polyhedron $C$,
a map $a:E\to C$, and
a homotopy equivalence $h:D\nu\cup_a C\to M$.

Two Poincar\'e embeddings $\P=(\nu,C,a,h)$ and $\P_1=(\nu_1,C_1,a_1,h_1)$ are called {\bf equivalent}, if there exist a bundle equivalence $b:\nu\to\nu_1$ and a homotopy equivalence $r:C\to C_1$ such that
$r\circ a\simeq a_1\circ b$ and $h\simeq h_1\circ h_{br}$, where the homotopy equivalence
$h_{br}:D\nu\cup_a C\to D\nu_1\cup_{a_1}C_1$ is constructed from $b$ and $r$ in the obvious way.

The Poincar\'e embeddings of isotopic embeddings are equivalent.
So the equivalence class of the Poincar\'e embedding is an isotopy invariant.

\proclaim{Browder-Wall Conjecture 9}
Suppose that $N$ and $M$ are closed smooth $n$- and $(n+k)$-manifolds, where $k\ge3$.

(a) If $n\ge 2k+3$ and $\P$ is a Poincar\'e embedding of $N$ in $M$, then there exists a smooth embedding $f:N\to M$ such that $\P(f)\sim\P$.

(b) If $n\ge 2k+4$ and $f,g:N\to M$ are smooth embeddings such that $\P(f)\sim\P(g)$, then $f$ and $g$ are smoothly isotopic.
\endproclaim

Although the conjecture in this explicit form was not published by Browder and Wall, we attribute it to them because it is very much in the spirit of [Br68, Theorem 2] [Wa70, \S11].
We believe the conjecture can be proved by invoking arguments from [Br68, Theorem 2], [Wa70, \S11], [Ma80, \S10].

For analogues in PL locally flat and TOP locally flat categories see [Wa70].

\head Proof of Lemma 7 \endhead

Since $\Sigma\Sy\sim\Sy(F)$, there is a bundle isomorphism $b:\nu\oplus1\to\nu(F)$ and a homotopy equivalence
$r:\Sigma C\to C(F)$ such that $r\circ a'\simeq a(F)\circ b$.

First, we construct {\it an $n$-manifold $N'\subset S^{n+k}$ and a map
$g:N'\to N$ of degree 1 such that $\nu':=\nu_{S^{n+k}}(N')=g^*\nu$.}
Since $\Sigma\Sy\sim\Sy(F)$, we have 
$$S^{n+k+1}\ \simeq\ D(\nu\oplus1)\cup_{a'}\Sigma C\ \simeq
\ \Sigma D\nu\cup_{\Sigma a}\Sigma C\ \simeq
\ \Sigma(D\nu\cup_aC).$$
Hence $H_i(D\nu\cup_aC;\Z)\cong H_i(S^{n+k};\Z)$. 
Since $D\nu$, $C$ and $S\nu$ are simply-connected,
using Seifert-van Kampen theorem one can show that $D\nu\cup_aC$ is simply-connected.
Therefore $D\nu\cup_aC\simeq S^{n+k}$.

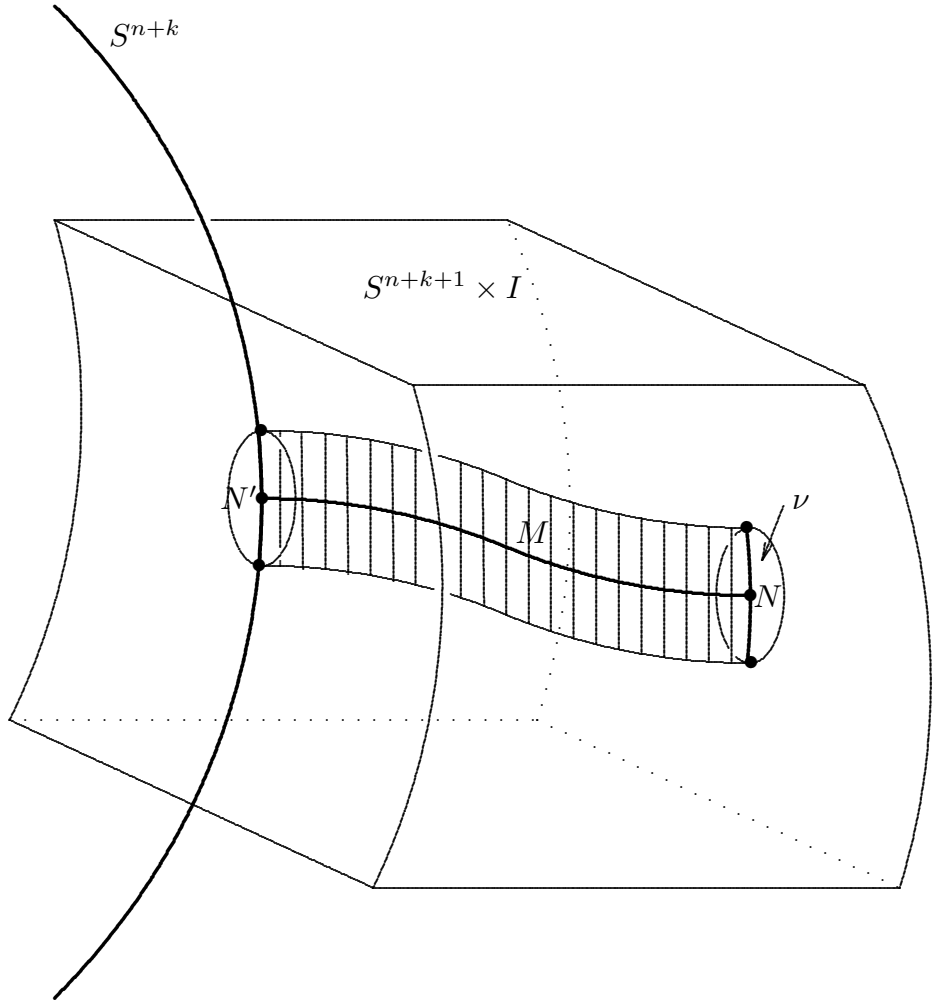
\input pictex
\beginpicture 
\setcoordinatesystem units <1.000mm,1.000mm>
\setplotarea x from -80 to 80, y from -80 to 80
\setsolid
\setplotsymbol ({\fiverm.})
\circulararc 15.806 degrees from -8.024 -4.728 center at -27.135 -73.161
\circulararc 5.263 degrees from 0.000 -7.500 center at -27.060 -72.828
\circulararc 14.854 degrees from -9.420 10.647 center at -26.976 -57.436
\circulararc 6.434 degrees from -0.000 7.500 center at -27.060 -57.828
\plot -32.637 36.225 0.167 36.225 /
\plot -50.000 36.225 -34.637 36.225 /
\circulararc 42.141 degrees from -55.000 -19.203 center at -124.434 15.000
\circulararc 42.309 degrees from -14.722 -37.813 center at -84.460 -4.235
\circulararc 42.276 degrees from 43.622 -37.813 center at -30.340 -15.000
\plot -14.722 -37.813 43.622 -37.813 /
\plot -55.000 -19.203 -14.722 -37.813 /
\plot -10.284 17.875 39.732 17.875 /
\plot 0.167 36.225 39.732 17.875 /
\plot -50.000 36.225 -10.284 17.875 /
\arrow <2mm> [.2,.4] from 30.724 4.497 to 28.376 -1.499
\plot 20.000 2.471 20.000 -12.529 /
\plot 15.000 3.153 15.000 -11.847 /
\plot 10.000 4.206 10.000 -10.794 /
\plot 5.000 5.647 5.000 -9.353 /
\plot 0.000 7.500 0.000 -7.500 /
\plot -5.000 9.353 -5.000 -5.647 /
\plot -10.000 10.794 -10.000 -4.206 /
\plot -15.000 11.847 -15.000 -3.153 /
\plot -20.000 12.529 -20.000 -2.471 /

\plot -2.591 8.453 -2.430 -6.546 /
\plot -12.621 11.233 -12.460 -3.766 /
\plot -17.630 12.045 -17.469 -2.954 /
\plot -22.635 12.494 -22.474 -2.505 /
\plot 22.500 -12.735 22.500 2.265 /
\plot 17.500 -12.233 17.500 2.767 /
\plot 12.500 -11.367 12.500 3.633 /
\plot 7.500 -10.123 7.500 4.877 /
\plot 2.500 -8.480 2.500 6.520 /

\ellipticalarc axes ratio 1:2 360 degrees from -23.310 5.383 center at -27.060 5.383
\ellipticalarc axes ratio 1:2 232 degrees from 25.500 -12.203 center at 27.060 -5.383
\ellipticalarc axes ratio 1:2 98 degrees from 24.500 0.097 center at 27.060 -5.383
\circulararc 22.500 degrees from 0.000 -7.500 center at 27.060 57.828
\circulararc 22.500 degrees from 0.000 7.500 center at 27.060 72.828
\plot -25.000 -2.147 -25.000 -1.600 /
\plot -25.000 -0.047 -25.000 10.753 /
\plot -25.000 12.306 -25.000 12.853 /
\plot 25.000 2.147 25.000 -12.853 /
\setplotsymbol ({\rm.})
\circulararc 22.500 degrees from 0.000 0.000 center at -27.060 -65.328
\circulararc 22.500 degrees from 0.000 0.000 center at 27.060 65.328
\circulararc 90.566 degrees from -50.000 -50.000 center at -104.460 5.000
\circulararc 11.108 degrees from -27.386 -2.110 center at -104.460 5.000
\circulararc 11.108 degrees from 26.733 -12.875 center at -50.340 -5.765
\setdots
\setplotsymbol ({\fiverm.})
\circulararc 42.039 degrees from 3.453 -19.203 center at -70.314 4.235
\plot -55.000 -19.203 3.453 -19.203 /
\plot 3.453 -19.203 43.622 -37.813 /
\put {$S^{n+k}$} at -40.000 57.000
\put {$S^{n+k+1} \times I$} at -7.045 28.951
\put {$M$} at 2.935 1.642
\put {$N$} at 29.000 -5.500
\put {$N'\ $} at -29.000 5.500
\put {$\nu$} at 32.681 4.925
\put {$\bullet$} at -27.386 -2.110
\put {$\bullet$} at -27.060 5.383
\put {$\bullet$} at -27.140 12.883
\put {$\bullet$} at 27.140 -12.883
\put {$\bullet$} at 27.060 -5.383
\put {$\bullet$} at 26.659 2.106
\endpicture


Let $g:S^{n+k}\to D\nu\cup_a C$ be a homotopy equivalence that is transverse regular on $N$.
Then $N'=g^{-1}(N)\subset S^{n+k}$ and $g|_{N'}$ are the required manifold and map, respectively.

Second, we construct a cobordism
$$M\subset S^{n+k+1}\times I\quad\text{between}
\quad N\subset S^{n+k+1}\times 0\quad\text{and}
\quad N'\subset S^{n+k}\times1 \subset S^{n+k+1}\times1$$
and a map
$$G:M\to N\quad\text{such that}\quad G|_N=\text{id},
\quad G|_{N'}=g\quad\text{and}
\quad\mu:=\nu_{S^{n+k+1}\times I}(M)=G^*\nu\oplus1.$$
$$\CD S^{n+k+1} @>>\Sigma g> \Sigma(D\nu\cup_aC) \\
@VV \eta' V @VV c V \\
T(\nu'\oplus1) @>> (g_*\oplus1)_T > T(\nu\oplus1) @>> b_T > T\nu(F)
\endCD$$
In order to construct such a cobordism consider the diagram above, where

$\bullet$ $g:S^{n+k}\to D\nu\cup_a C$ is the homotopy equivalence constructed using the first step,

$\bullet$  $\eta'$ is the collapsing map (representing a normal invariant of $N'$), and

$\bullet$ $c$ is the collapsing of $\Sigma C$ to a point.

Clearly, the diagram is commutative.
One can check that $b_T\circ c\circ\Sigma g\sim\alpha(F)$.
Hence there is a homotopy $G:S^{n+k+1}\times I\to T\nu(F)$, transverse regular on $N$, between
$G_0=b_T\circ(g_*\oplus1)_T\circ\eta'$ and a map $G_1$ representing `the normal invariant' $\alpha(F)$.
Therefore $G_1^{-1}(N)=N$ and the restriction
$G_1|_N=$id is covered by a bundle isomorphism of $\nu(F)$ (which we assume to be the identity).
So $M:=G^{-1}(N)$ is a $(n+1)$-dimensional cobordism between $G_0^{-1}(N)=N'$
and $G_1^{-1}(N)=N$.
From now on, denote $G|_M:M\to N$ by $G$.
Clearly, $M$ and $G$ are as required in the second step.

\proclaim{Claim} There exist a manifold $X$, with boundary $\partial X=S\nu$,
that is cobordant to
$$X_1:=(S^{n+k}-\delet D\nu')\bigcup\limits_{S\nu'=S(\mu|_{N'})}S\mu$$
modulo the boundary, and a homotopy equivalence $R':X\to C$ whose restriction
to $\partial X=S\nu$ is $a$.
\endproclaim

\demo{Proof of Lemma 7 modulo Claim}
Let
$$W:=X\bigcup\limits_{S(\mu|_N)=S\nu}D\nu\quad\text{and}
\quad W_1:=X_1\bigcup\limits_{S(\mu|_N)=S\nu}D\nu\ =
\ \partial(D^{n+k+1}\cup D\mu)$$
Observe that $W$ contains $N$ with the normal bundle $\nu$.
Since $W$ is cobordant to $W_1$, it follows that $W$ is null-cobordant.
Also, $W$ is mapped to $D\nu\cup_a C \simeq S^{n+k}$ by a homotopy equivalence
(which is the `union' of $R'$ and a diffeomorphism on $D\nu$).
The stable normal bundle of $W$ is trivial (because it is induced from the
pull-back of the trivial normal bundle of $S^{n+k}$ in $S^{n+k+1}$).
The normal invariant of $W$ is also trivial.
Hence by the Novikov surgery theorem $W$ is a homotopy sphere and is the boundary of a parallelizable manifold.  Therefore the connected sum $W'\#\Sigma$ with the inverse homotopy sphere (away from $D\nu$) is
diffeomorphic to $S^{n+k}$.
We have $N\subset D\nu\subset W\#\Sigma$ and the normal system of $N$ in $W\#\Sigma$ is $(\nu,C,a)$.
\qed\enddemo

\demo{Proof of Claim}
We may assume that $a:S\nu\to C$ is an inclusion.
Define a map
$$R:D\mu\cup S^{n+k}\times I\to D\nu\cup_a C\qquad\text{by}\qquad
\cases R(x,s)=g(x) & (x,s)\in S^{n+k}\times I\\
R(y)=G_*(y) & y\in D\mu\endcases.$$
We have
$$R^{-1}(N)=(N'\times I)\cup M,\quad\text{so} \quad R^{-1}(N)\cap
W=N\subset M,\quad R|_{M}=G \quad\text{and}\quad R|_N=\text{id}.$$
So the restriction of the above $R$ is a new map $R:X_1\to C$
extending $a$ from $\partial X_1=S\nu$.
The space $S(\mu\oplus1)$ can be identified with the boundary of a regular
neighbourhood of $M$ in $S^{n+k+1}\times I$.
Hence $S\mu\subset S(\mu\oplus1)$ is embedded into $S^{n+k+1}\times I$.
Since also $S^{n+k}-\delet D\nu'$ is embedded in $S^{n+k+1}$ (with trivial
normal bundle), it follows that the manifold $X_1$ can be embedded in
$S^{n+k+1}\times I\subset S^{n+k+2}$.
We have $\nu_{S^{n+k+1}\times I}(S\mu)=\nu_{S(\mu\oplus1)}(S\mu)\oplus G^*1$ is trivial.
Therefore the stable normal bundle of $X_1$ is also trivial.

Let $W_1=X_1\bigcup\limits_{S(\mu|_N)=S\nu}D\nu$.
Clearly, $W_1$ contains $N$ with normal bundle $\nu$.
Recall that $W_1$ is null-cobordant.
Since $W_1$ satisfies the Poincar\'e duality,
$D\nu\cap C=S\nu$ satisfies the Poincar\'e duality and
$H_{n+k-1}(D\nu)=0=H_{n+k-1}(C)$, it follows that $(C,S\nu)$ is a
Poincar\'e pair in dimension $n+k$.
If $n+k=\dim X_1$ is odd, there is no obstruction to modify $R$ by surgery
to a map $R':X\to C$ which is a homotopy equivalence extending $a$.
If $n+k=2l$, then  by making surgery we may assume that $R_*=R(X_1)_*$
induces an isomorphism in homology up to dimension $l-1$ and an
epimorphism in dimension $l$.
There is however a surgery obstruction $\sigma(R(X_1))$ to making $R(X_1)_*$ a monomorphism in dimension $l$.
There is analogous obstruction $\sigma(R(W_1))$ to surgery of
$R(W_1):W_1\to S^{2l}$ to a homotopy equivalence with $D\nu\cup_aC$.
We have the following commutative diagram whose vertical arrows are induced by inclusion:
$$\CD \Ker R(X_1)_* @>> > H_l(X_1) @>> R(X_1)_* > H_l(C) @>> > 0\\
@VV \cong V @VV V @VV V\\
\Ker R(W_1)_* @>> = > H_l(W_1) @>> R(W_1)_* > H_l(S^{2l}) @>> > 0 \endCD$$
Hence by definition of $\sigma$ [Wa70] we have $\sigma(R(X_1))=\sigma(R(W_1))$.
However, $\sigma(R(W_1))=0$ since $W_1$ is cobordant to $S^{n+k}\simeq D\nu\cup_aC$.
Since $R|_{\partial X_1}$ is a diffeomorphism, we can do surgery in the
interior of $X_1$ and leave $\partial X_1$ and $R|_{\partial X_1}$
unchanged.
\qed\enddemo

{\bf References}

[Al24] J.W. Alexander, On the Subdivision of 3-Space by Polyhedron, Proc. Natl. Acad. Sci. USA 10, 6-8 (1924).

[Al24'] J.W. Alexander, An Example of a Simply Connected Surface Bounding a Region Which Is Not Simply Connected, Proc. Natl. Acad. Sci. USA 10, 8-10 (1924).

[ARS] P. Akhmetiev, D. Repov\v s, and A. Skopenkov, Embedding Products of Low-Dimensional Manifolds in $\R^m$,
Topology Appl. 113, 7-12 (2001).

[Ba65] D. Barden, Simply Connected Five-Manifolds, Ann. Math., Ser. 2, 82, 365-385 (1965).

[Br60] M. Brown, A Proof of the Generalized Schoenflies Theorem, Bull. Am. Math. Soc. 66, 74-76 (1960).

[Br68] W. Browder, Embedding smooth manifolds, in Proc. Int. Congr. Math., Moscow, 1966 (Mir, Moscow, 1968),
pp. 712-719.

[CS08] D. Crowley and A. Skopenkov. A classification of smooth embeddings of 4-manifolds in 7-space, II, Intern. J. Math., 22:6, 731-757 (2011), arxiv:math/0808.1795.

[CS16I] D. Crowley and A. Skopenkov, Embeddings of non-simply-connected 4-manifolds in 7-space. I. Classification modulo knots, Moscow Math. J., 21 (2021). arXiv:1611.04738.

[CS16II] D. Crowley and A. Skopenkov, Embeddings of non-simply-connected 4-manifolds in 7-space. II. On the smooth classification, Proc. A of the Royal Soc. of Edinburgh, to appear. arXiv:1612.04776.

[DV09] R.J. Daverman and G.A. Venema, Embeddings in Manifolds, GSM 106, AMS,  Providence, RI, 2009.

[FF89] A.T. Fomenko and D.B. Fuchs, Homotopical Topology, Springer, Berlin, 2016.

[Gi71] S. Gitler, Immersion and Embedding of Manifolds, Proc. Symp. Pure Math. 22, 87-96 (1971).

[Gl62] H. Gluck, The Embedding of Two-Spheres in the Four-Sphere, Trans. Am. Math. Soc. 104 (2), 308-333 (1962).

[GS99] R. Gompf and A. Stipsicz, 4-Manifolds and Kirby Calculus, GSM 20, AMS, Providence, RI, 1999.

[Ha37] W. Hantzsche, Einlagerung von Mannigfaltigkeiten in euklidische R\" aume, Math. Zeitschrift 43:1, 38-58 (1937).

[Ha66] A.~Haefliger, Enlacements de spheres en codimension superiure a 2, Comment. Math. Helv. 41, 51-72 (1966-67).

[Ha68] A. Haefliger, Knotted Spheres and Related Geometric Topic, in Proc. Int. Congr. Math., Moscow, 1966 (Mir,
Moscow, 1968), pp. 437-445.

[Ha86] N. Habegger, Knots and links in codimension greater than 2, Topology, 25:3,  253-260 (1986).

[Hi59] M. W. Hirsch, Immersions of Manifolds, Trans. Am. Math. Soc. 93, 242-276 (1959).

[Ho40] H. Hopf, Systeme symmetrischen Bilinearformen und euklidische Modelle der projectiven R\"aume, Vierteljahresschr. Naturforsch. Ges. Z¨urich 85, 165-178 (1940).


[IF] http://www.map.mpim-bonn.mpg.de/Intersection\_form

[Im] http://www.map.mpim-bonn.mpg.de/Immersion\#The\_smooth\_case

[Is] http://www.map.mpim-bonn.mpg.de/Isotopy

[Le63] J. Levine, On Differentiable Embeddings of Simply-Connected Manifolds, Bull. Am. Math. Soc. 69, 806-809
(1963).

[Le65] J. Levine, Unknotting Spheres in Codimension 2, Topology 4, 9-16 (1965).

[LF] http://www.map.mpim-bonn.mpg.de/Linking\_form

[LS69] W.B.R. Lickorish and L.C. Siebenmann, Regular Neighborhoods and the Stable Range, Trans. Am. Math.
Soc. 139, 207-230 (1969).

[Ma59] B. Mazur, On Embeddings of Spheres, Bull. Am. Math. Soc. 65, 91-94 (1959).

[Ma80] R. Mandelbaum, Four-Dimensional Topology: An Introduction, Bull. Am. Math. Soc. 2, 1-159 (1980).

[Mo60] M. Morse, A Reduction of the Schoenflies Extension Problem, Bull. Am. Math. Soc. 66, 113-117 (1960).

[No64] S.P. Novikov, Homotopy Equivalent Smooth Manifolds, Izv. Akad. Nauk SSSR, Ser. Mat. 28, 365-474 (1964).

[Pa57] C.D. Papakyriakopoulos, On Dehn's Lemma and the Asphericity of Knots, Ann. Math., Ser. 2, 66, 1-26
(1957).

[RBS] D. Repov\v s, N. Brodsky, and A.B. Skopenkov, A Classification of 3-Thickenings of 2-Polyhedra, Topology Appl. 94, 307-314 (1999).

[RS72] C.P. Rourke and B.J. Sanderson, Introduction to Piecewise-Linear Topology (Springer, Berlin, 1972), Ergebn. Math. 69.

[RS99] D. Repov\v s and A.B. Skopenkov, New results on embeddings of polyhedra and manifolds into Euclidean spaces, Russ. Math. Surv. 54:6, 1149-1196 (1999).


[SE62] N.E. Steenrod and D.B.A. Epstein, Cohomology Operations, Princeton Univ. Press, Princeton, 1962.

[Sk05] A.  Skopenkov, A classification of smooth embeddings of 4-manifolds in 7-space, I, Topol. Appl., 157,   2094-2110 (2010). arXiv:math/0512594.

[Sk06] A. Skopenkov, Embedding and knotting of manifolds in Euclidean spaces, London Math. Soc. Lect. Notes, 347,  248-342 (2008). arXiv:math/0604045.

[Sk06'] A. Skopenkov, A classification of smooth embeddings of 3-manifolds in 6-space, Math. Zeitschrift, 260:3,  647-672 (2008). arxiv:math/0603429.

[Sk16] A. Skopenkov,  Embeddings in Euclidean space: an introduction to their classification, to appear in Boll. Man. Atl. http://www.map.mpim-bonn.mpg.de/
\linebreak
Embeddings\_in\_Euclidean\_space:\_an\_introduction\_to\_their\_classification

[Sk20] A. Skopenkov, Algebraic Topology From Geometric Viewpoint (in Russian), MCCME, Moscow, 2020 (2nd edition). Electronic version:

http://www.mccme.ru/circles/oim/home/combtop13.htm\#photo. Accepted for the English translation by Springer.

[Sk20u] A. Skopenkov, A user's guide to basic knot and link theory, in: `Topology, Geometry, and Dynamics:
Rokhlin Memorial'. Contemp. Math. 772,  AMS, Providence, RI, 2021.
Russian version to appear in Mat. Prosveschenie. arXiv:2001.01472.

[Sm61] S. Smale, Generalized Poincar\'e's Conjecture in Dimensions Greater than 4, Ann. Math., Ser. 2, 74, 391-466
(1961).

[St63] J. Stallings, On Topologically Unknotted Spheres, Ann. Math., Ser. 2, 77, 490-503 (1963).

[SW] http://www.map.mpim-bonn.mpg.de/Stiefel-Whitney\_characteristic\_classes

[Th51] R. Thom, Une th\'eorie intrins\`eque des puissances de Steenrod, in Colloq. Topol., Strasbourg, 1951, 
no. 6.

[Wa66] C. T. C. Wall, Classification Problems in Differential Topology. IV: Thickenings, Topology 5, 73-94 (1966).

[Wa70] C. T. C. Wall, Surgery on Compact Manifolds (Academic, London, 1970; AMS, Providence, RI, 1998).

[Ze93] E. C. Zeeman, A Brief History of Topology, UC Berkeley, October 27, 1993, On the occasion of Moe Hirsch's 60th birthday,

http://zakuski.utsa.edu/\~{}gokhman/ecz/hirsch60.pdf.

\enddocument